\documentclass{amsart}
\usepackage{hyperref,amsmath,amsthm,amsfonts,amscd,flafter,epsf}
\begin{document}

\title{Heegaard Floer Homologies of Pretzel Knots}
\author{Eaman Eftekhary}
\begin{abstract}
We compute the knot invariants defined by Ozsv\'{a}th and Szab\'{o} \cite{OS-knot} for $(-2a,2b+1,2c+1)$ and $(2a,-(2b+1),2c+1)$-Pretzel knots, for
$a,b,c \in \mathbb{Z}^+$.
We will also deduce certain information on the contribution of homotopy disks to the boundary counts.
\end{abstract}
\maketitle
\section{Introduction}
Recently, Ozsv\'{a}th and Szab\'{o} have introduced several Floer homology constructions, assigning invariants to $3,4$-manifolds and more recently to knots,
links and contact structures on $3$-manifolds (\cite{OS-3m1,OS-3m2,OS-4m,OS-knot,OS-fibered}
, see also \cite{Ras2}).\\

The common flavor of these constructions is assigning a surface $\Sigma_g$ of genus $g$ and two sets of disjoint, homologically independent circles
$(\alpha_1,\alpha_2,...,\alpha_g),(\beta_1,\beta_2,...,\beta_g)$ on $\Sigma_g$ together with some marked points in the different regions of
$\Sigma_g-\cup_{i=1}^{g}\alpha_i -\cup_{i=1}^g \beta_i$ to the object in hand, which may be the $3$-manifold, the knot or the contact structure.
Then they consider the Lagrangian Floer homology associated with the submanifolds $\mathbb{T}_\alpha=\alpha_1\times ...\alpha_g$ and
$\mathbb{T}_\beta=\beta_1\times ...\times \beta_g$ in the symmetric product:

\begin{equation}
Sym^g(\Sigma_g)=(\prod_{i=1}^g \Sigma_g) /S_g
\end{equation}
where $S_g$ is the symmetric group on $g$ letters. The associated Floer homologies are then proved to be independent of the various choices made meanwhile.\\

In particular the knot Floer homology, called the Heegaard Floer homology of knots by Ozsv\'{a}th and Szab\'{o}, is easy to compute at least up to the level of
chain complexes. In fact, given a knot projection diagram, one can fix a marked point on the knot diagram and consider the associated Kauffman states
(see \cite{OS-knot} for definitions) of the knot projection. These Kauffman states will be unions of marked points, one in each of the regions formed
by the
knot projection (except for the neighbors of the marked point on the knot). The marked points are put close to the self intersections of the diagram of knot
projection. We assume for each self intersection that exactly one of the $4$ neighboring regions has its marked point close to the chosen
intersection. The chain complex will be freely generated by these Kauffman states.
\\

In the standard knot theory, one will assign two types of gradings to each Kauffman state, which Ozsv\'{a}th and Szab\'{o} call the Maslov and spin grading
respectively, because of the relation to the $3$-manifold invariants and the notation used there (\cite{OS-knot}). These can be read off easily (combinatorially)
from the diagram. The boundary maps will respect the spin grading and will lower the Maslov grading by one. So there will be a double grading imposed on the
Floer homology groups. We will denote by $\widehat{\text{HFK}}(K,s)$ the Heegaard Floer homology of the knot $K$ in the spin grading $s$ and by
$[\widehat{\text{HFK}}(K,s)]_m$ the part of these homology groups in bi-grade $(m,s)$, where $m,s$ denote the Maslov and spin gradings, respectively.\\

Computation of the boundary maps is harder though. By thickening the knot projection, one will get a surface. Associated with the domains in the
knot projection, the $\beta$-curves are drawn on this surface. The $\alpha$-circles will be associated with the self intersections and one
special  $\alpha$-circle will be assigned to the marked point. The counts for the boundary operator will reduce to deciding for certain combinations
$\mathcal{D}=\sum_{i}\lambda_i.D_i$ of the regions $D_i$ of $\Sigma_g-\cup_{i=1}^{g}\alpha_i -\cup_{i=1}^g \beta_i$, that how many holomorphic maps
$u:\Delta \rightarrow Sym^g(\Sigma_g)$ have a domain equal to $\mathcal{D}$. Here $\Delta$ denotes the unit disk (see \cite{OS-knot} for the full
description).\\

In general certain domains are discovered to support holomorphic disks and their contributions to the boundary operators have been computed. A class of
them is illustrated in \textbf{figure.7}. There are also long exact sequences assigned to a triple $(L_{+},L_{-},L_{0})$ where $L_-$ is obtained by
changing a positive intersection in a projection diagram of $L_+$ to a negative one, and $L_0$ is obtained by resolving this intersection so that the
orientations are preserved in different pieces of the resulting link. This exact sequence may be used to get information on the contribution of the topological
disks to the boundary operator or help to compute the homology in some cases. We refer to \cite{OS-knot} for the precise statement, since we will not be using
this exact sequence.\\

Using these exact sequences and the information about the holomorphic representatives of certain domains, several computations have been done (c.f. \cite{Ras,Ras2,OS-knot}).
It turned out that the Heegaard Floer homology is completely determined  knowing the
symmetrized Alexander polynomial $\Delta_K(t)$ of an alternating knot $K$, together with the signature $\sigma(K)$. Moreover the homology is
supported on a line in the $(s,m)$-plane.
There are also computations for the torus knots and the $(2m+1,2n+1,2k+1)$-pretzel knots ($m,n,k\in \mathbb{Z}$) and some partial but very useful
information on the fibered knots. \\

In this paper we deal with the remaining three stranded pretzel knots.
Note that the only remaining cases are the $(-2a,2b+1,2c+1)$ and $(2a,-2b-1,2c+1)$ for
$a,b,c \in \mathbb{Z}^+$. The rest of them are either the mirrors of these, or lie in the category of previous computations.\\

The result is that again the Floer homologies are as simple as possible. Denote the $(m,n,k)$-pretzel knot by $K(m,n,k)$ then part of the result is the
following:\\

\textbf{Theorem.1.} \emph{Suppose that $a,b,c \in \mathbb{Z}^+, a\leq b\leq c$. Let $K=K(-2a,2b+1,2c+1)$. Then writing the spin grading $s$ as $s=b+c-p=b-c-q-1$
we will have:
\begin{displaymath}
1)\ \ \ [\widehat{\text{HFK}}(K,b+c+1)]_0= [\widehat{\text{HFK}}(K,-(b+c+1))]_{-2(b+c+1)}  =\oplus_{i=1}^a \mathbb{Z}
\end{displaymath}
Also we will have:
\begin{displaymath}
\begin{array}{ccc}
2)  & [\widehat{\text{HFK}}(K,s)]_{-(p+1)}=\oplus_{i=1}^{(2a-1)-p} \mathbb{Z} \ \ \  & 0\leq p\leq 2a-1\\
3)  & [\widehat{\text{HFK}}(K,s)]_{-p}=\oplus_{i=1}^{p-(2a-1)} \mathbb{Z}      & 2a\leq p \leq 2b \\
4)  & [\widehat{\text{HFK}}(K,s)]_{-p}=\oplus_{i=1}^{2(b-a)+1} \mathbb{Z}      & 2b <p \leq 2c \\
5)  & [\widehat{\text{HFK}}(K,s)]_{-2c-p+1}=\oplus_{i=1}^{2(b-a)-p} \mathbb{Z} & 0\leq p\leq 2(b-a) \\
6)  & [\widehat{\text{HFK}}(K,s)]_{-2c-p}=\oplus_{i=1}^{p-2(b-a)} \mathbb{Z}   & 2(b-a)< p< 2b\\
\end{array}
\end{displaymath}
and the rest of homology groups will be zero.}\\
\\
Note that this is completely determined by stating the following: ``The Floer homology is supported on the two lines $s=m+(b+c+1)$ and $s=m+(b+c)$. For each $s$
only on one of the two lines we may have a nontrivial group depending on the sign of $(-1)^s a_s$. This nontrivial group
will be a sum of $|a_s|$ copies of $\mathbb{Z}$, where $\Delta_K(t)=\sum_{s}a_s t^s$ is the symmetrized Alexander polynomial.''. This statement remains true even without the assumption
$a\leq b \leq c$.\\

In a similar way we state
the next theorem as follows:\\

\textbf{Theorem.2} \emph{For $K=K(2a,-2b-1,2c+1)$, with $a,b,c \in \mathbb{Z}^+$ let $\Delta_K(t)=\sum_{s}a_s t^s$ be the symmetrized Alexander polynomial
of $K$. Then the Heegaard Floer homology is supported on the two lines $s=m+(c-b)-1,s=m+(c-b)$. For each $s$ on one of the two lines we will
have $\oplus_{i=1}^{|a_s|}\mathbb{Z}$, depending on  the sign of $(-1)^s a_s$. The rest of Heegaard Floer homology groups will be zero.}\\
\\
Note that this completely determines the Floer homology of the knot $K$.\\

\begin{figure}
   \mbox{\vbox{\epsfbox{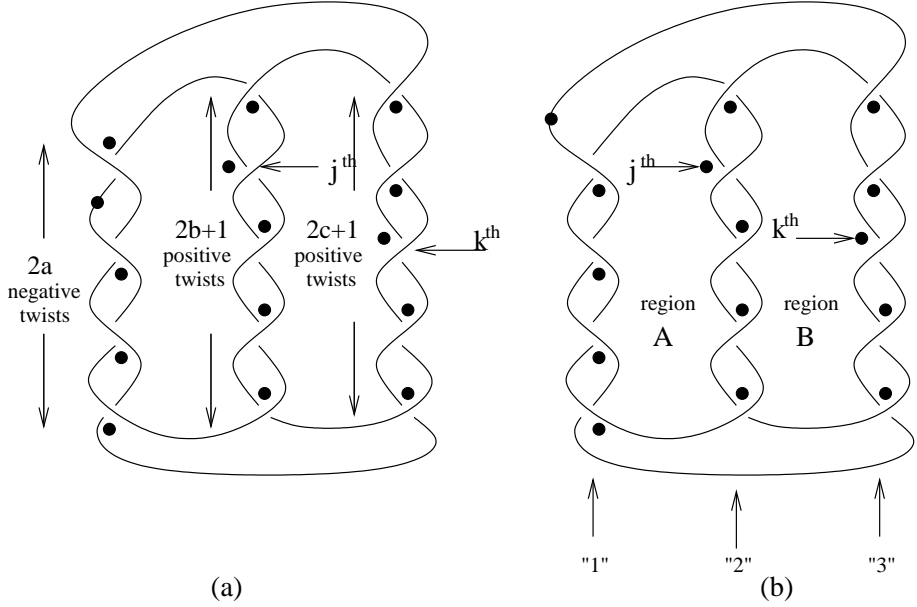}}}
   \caption{\label{fig:figure1}
   {A $(-2a,2b+1,2c+1)$-pretzel knot for $a=b=c=2$. Note that the twists on the first pair of strands are negative twists. The Kauffman states that
    are shown are $A_{jk}(a)$ and $A_{jk}(b)$.}}
\end{figure}

\textbf{Acknowledgment.} I would like to thank Peter Ozsv\'{a}th and Zolt\'{a}n Szab\'{o} for many helpful discussions.

\section{Heegaard Diagram and Kauffman States}

We remind the reader that a $(m,n,k)$-pretzel knot is a knot obtained from 3 pairs of strands which are twisted $m,n$ and $k$ times respectively
and then the ends are glued together cyclically as is shown in \textbf{figure.1} (also see \cite{knot}).\\

We will only discuss the proof in the $(-2a,2b+1,2c+1)$ case with $a\leq b \leq c$. The other cases are completely similar.\\

In order to do the computations, first one should fix a pointed Heegard diagram of the knot. The diagram that we will be using is the standard one, obtained
from the plane projection of the pretzel knot. Namely, we thicken the projection of the $(-2a,2b+1,2c+1)$-pretzel knot to obtain a genus $2(a+b+c)+3$
surface.\\
Then the $\beta$-circles correspond to different finite regions in the knot projection. They go around the ``whole'' corresponding to these regions in the
thickened knot once. The $\alpha$-circles correspond to the intersection points in the knot projection and cut all of the $\beta$-circles that are associated
with the neighboring areas of the intersection, once. \textbf{Figure.2} shows the construction for the trefoil knot .
We refer the reader to \cite{OS-knot} for a careful description of this Heegaard diagram for arbitrary knots.\\

\begin{figure}
   \mbox{\vbox{\epsfbox{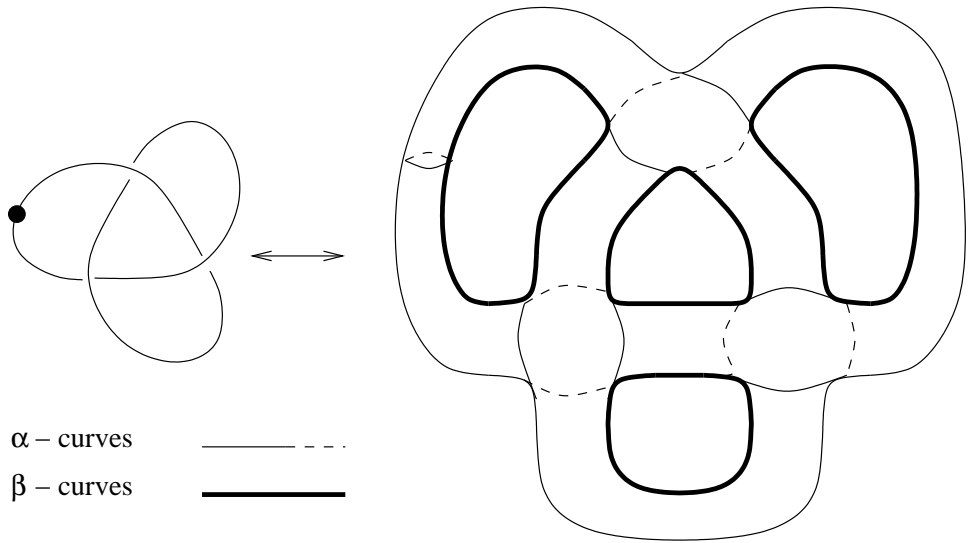}}}
   \caption{\label{fig:figure1}
   {To each knot projection is associated a Heegaard diagram. A sample Heegaard diagram for the trefoil with a marked point is illustrated. The bold lines
denote the $\beta$ curves and the small loops are $\alpha$-curves. There is a special $\alpha$-curve associated to the marked point.}}
\end{figure}

There is a special choice for an additional (exceptional) $\alpha$-curve which will go around one of the handles in the knot projection.
 This is in correspondence with a choice of a marked point somewhere ``on'' one of the arcs of the knot projection. Here we are assuming that one of the
two neighboring regions of this arc in the knot projection is the unbounded region. The Floer homology is independent of the choice of this point.
We will choose two
different points, as suggested in parts (a) and (b) of \textbf{figure.1} and will study the homology groups formed out of the associated chain complexes
of these two diagrams. A comparison of the boundary maps in the two complexes will give us the result.\\

Below we describe the intersection points of $\mathbb{T}_\alpha$ and $\mathbb{T}_\beta$ (the associated tori of the $\alpha$ and $\beta$ circles in the
symmetric product of the surface). These are in correspondence with the Kauffman states of the knot, as described in \cite{OS-knot}.\\

A Kauffman state is a union of $g-1$ marked points (with $g$ the number of $\alpha, \beta$ curves) in the diagram of the knot projection. Each of the
finite regions will get a marked point except for the neighboring region of the special point on the knot. The marked points are put close to the self
intersections in the knot projection, i.e. in one of the four quadrants. We require that for each self intersection, exactly one of the $4$ quadrants
gets a marked point.\\

In \textbf{figure.1} we have marked two of the regions by the letters A and B. In both cases (a) and (b), any Kauffman state
will have a marked point in  region A and one in  region B. \\
In our knot projection, there are three pairs of strands twisted $-2a,2b+1$ and $2c+1$ times respectively and then the ends are cyclically glued back
together. Let us call these pairs ``1'', ``2'' and ``3'' respectively. There are several intersections in the knot projection associated with each of these
pairs which we may number from top to the bottom in the picture, starting as $0,1,2,...$ . So we will be denoting the $j^\text{th}$ intersection from the
top on the second pair of strands by $2[j-1]$, etc. . \\

With this notation fixed, there are three possibilities for the marked points in regions A,B:\\
\\
(1) The marked point in ``A'' is located next to an intersection on the second pair of strands, say $2[j]$, and the marked point on ``B'' is next to
$3[k]$.\\
(2) ``A'' has a marked point next to $1[i]$ and ``B'' has one next to $3[k]$.\\
(3) The marked point on ``A'' is next to $1[i]$ and that of ``B'' is next to $2[j]$.\\
Note that $i \in \{0,1,...,2a-1\}, j\in \{0,1,...,2b\}, k\in \{0,1,...,2c\}$.\\

It is important to note that in both cases (a) and (b), the position of the marked points in the regions ``A'' and ``B'' determines the Kauffman state
completely. We will call the state described in (1), $A_{jk}(a)$ or $A_{jk}(b)$ depending on whether we are considering the case with the marked point as
in \textbf{figure 1(a)} or the case shown in \textbf{figure 1(b)}.
Similarly from the description (2) above we will get $B_{ik}(a),B_{ik}(b)$ and from (3) will get $C_{ij}(a)$ and $C_{ij}(b)$ (see \textbf{figure.3}).
Note that $B_{0,k}(a)$
and $C_{0,j}(a)$ are special and somehow different from what is shown in \textbf{figure.3}.

If a statement is true for both $A_{jk}(a)$ and $A_{jk}(b)$ we will simply state it for $A_{jk}$, etc..\\

\begin{figure}
   \mbox{\vbox{\epsfbox{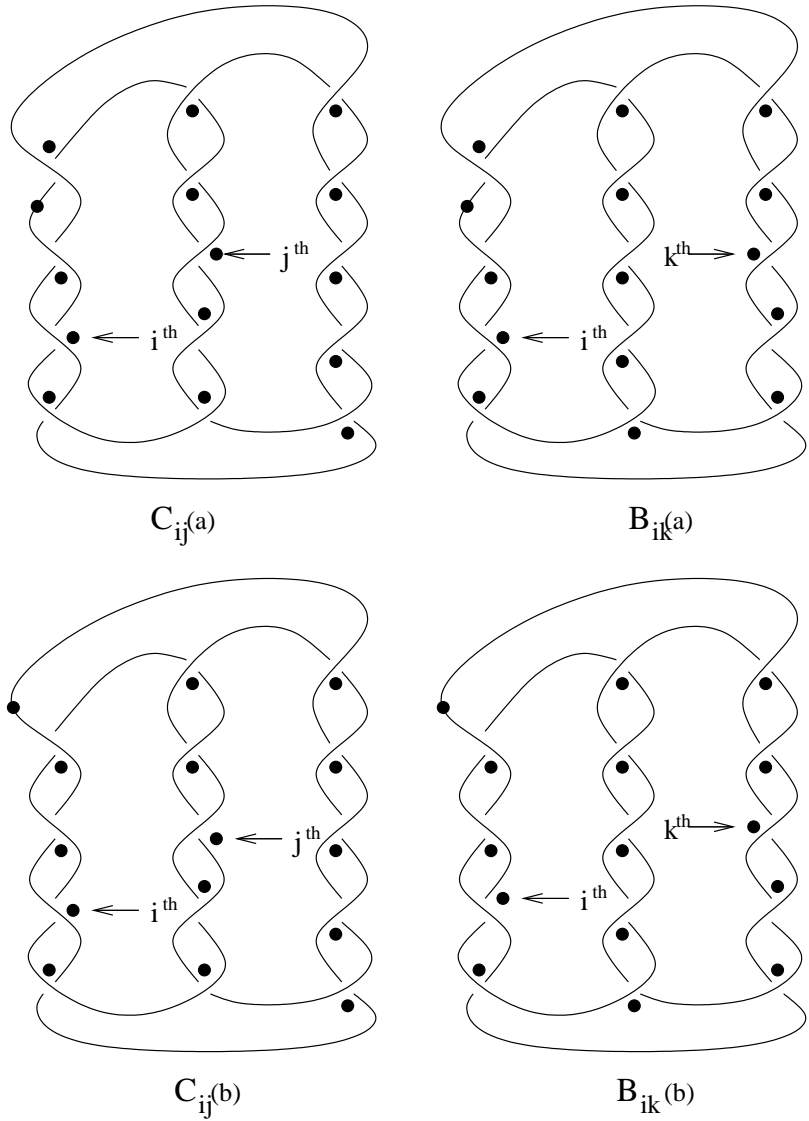}}}
   \caption{\label{fig:figure1}
   {Different types of Kauffman states are illustrated. The upper two are of type (a) and the next two are of type (b).}}
\end{figure}

\begin{figure}
   \mbox{\vbox{\epsfbox{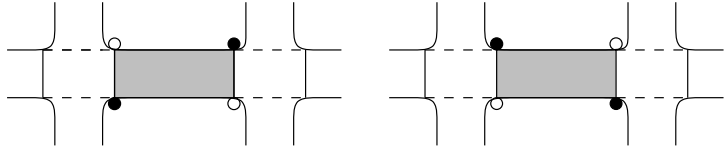}}}
   \caption{\label{fig:figure1}
   {Assigning small area to the shaded region, one may assume that the blank and bold Kauffman state are canceled against each other}}
\end{figure}

We now turn to compute the Maslov grading and the spin grading of these states. It is easy to see that all of the intersection points of
the knot projection are positive intersections and using the local contributions to Maslov and spin grading (as is shown in \textbf{figure.5})
we may compute them as follows (we will write $\epsilon(i)$ for $i-2[\frac{i}{2}]$):

\begin{equation}
  \left|\begin{array}{c}
   \hline\\
   \begin{array}{ccc}
      \begin{array}{cc}
          {\textbf{x}}& \\
            \hline \\
          A_{jk}(a)&  \\
          B_{ik}(a)&  i>0\\
          B_{0k}(a)&   \\
          C_{ij}(a)& i>0 \\
          C_{0j}(a)&
      \end{array}&
      \left|
      \begin{array}{c}
          m({\textbf{x}})\\
          \hline \\
          j-k-2b\\
          -\epsilon (i)-k\\
          -k-2b-1\\
          j-\epsilon (i)-2b-2c-1\\
          j-2b
      \end{array}
      \right|&
      \begin{array}{c}
          s({\textbf{x}})\\
          \hline \\
          (j-k)+(c-b)\\
          b+c+1-\epsilon(i)-k\\
           -k+(c-b)\\
          j-\epsilon(i)-b-c\\
          j-b+c+1
      \end{array}
   \end{array}\\
  \hline
 \end{array}\right|
\end{equation}
\\

\begin{equation}
  \left|\begin{array}{c}
   \hline\\
   \begin{array}{ccc}
      \begin{array}{c}
          {\textbf{x}} \\
            \hline\\
          A_{jk}(b)\\
          B_{ik}(b)\\
          C_{ij}(b)
      \end{array}&
      \left|
      \begin{array}{c}
          m({\textbf{x}})\\
          \hline \\
          j-k-2b\\
          -\epsilon (i)-k\\
          j-\epsilon (i)-2b-2c-1
      \end{array}
      \right|&
      \begin{array}{c}
          s({\textbf{x}})\\
          \hline \\
          (j-k)+(c-b)\\
          b+c+1-\epsilon(i)-k\\
          j-\epsilon(i)-b-c
      \end{array}
   \end{array}\\
  \hline
 \end{array}\right|
\end{equation}
\\

For an intersection point (or equivalently Kauffman state) \textbf{x}, define $\Delta(\textbf{x})=s(\textbf{x})-m(\textbf{x})-(b+c)$. Then from the
above table, we see immediately that $\Delta(A_{jk}(a))=\Delta(A_{jk}(b))=0$, while $\Delta(B_{ik})=\Delta(C_{ij})=1$. On the other hand, for two
Kauffman states \textbf{x,y}, we know that \textbf{y} may appear in the boundary of \textbf{x} only if $s(\textbf{x})=s(\textbf{y})$ and
$m(\textbf{x})=m(\textbf{y})+1$. If so, we will definitely have $\Delta(\textbf{y})-\Delta(\textbf{x})=1$. As a result the only boundary maps will
go from $A_{jk}$'s to $C_{ij}$'s and $B_{ik}$'s.

\begin{figure}
   \mbox{\vbox{\epsfbox{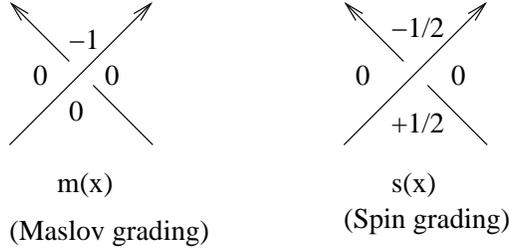}}}
   \caption{\label{fig:figure1}
   {Local contribution to Maslov and spin grading, for positive intersections.}}
\end{figure}

\section{Topological disks between the Kauffman states}

With the above computation in hand, the next step is to compute the coefficients of different regions in the associated domain of the topological
disks between the relevant Kauffman states. We will use the knot projections with pairs of integers on the edges, to denote the domain with these pairs
as the coefficients on the upper face and lower face of the corresponding handle of that edge (arc). For a more careful description of these
notations we refer the reader to \cite{OS-knot}.\\
There are three arcs on the top of our diagram of the knot projection which connect the ends of the pairs of strands to each other cyclically.
Correspondingly, there are two regions associated to each of these arcs, one on the upper face of the surface and one on the lower face.
We will compute all the possible topological disks between $A$-type and $B,C$-type Kauffman states with the property that their coefficient on the lower
face of the surface on these three arcs is zero. It is easy to see that
the only possible disks of Maslov index $1$ between $A_{jk}$'s and $B_{il}$'s in both cases are the ones going from $A_{2b-j,k}$ to
$B_{2a-j-1,j+k+\epsilon(j)}$ and the one going from $A_{0k}$ to $B_{0k}$.
This means that any other disk either has a negative coefficient in the region associated to the ``top arcs'', or it has at
least one strictly positive coefficient in one of the punctured ones.\\

We also note here that in case (b), the intersection points $B_{ik}(b)$ and $C_{ij}(b)$ live in different spin gradings, except for one case.
In fact, if
$s(B_{i_1k}(b))=s(C_{i_2j}(b))$ then:
\begin{equation}
b+c+1-\epsilon(i_1)-k=j-\epsilon(i_2)-b-c
\end{equation}
Or putting it in other form, it says that:
\begin{equation}
(2b-j)+(2c-k)+(1-\epsilon(i_1))+\epsilon(i_2)=0
\end{equation}
This can be the case only if $j=2b,k=2c$ and $i_1,i_2$ are respectively odd and even. In this case the spin grading will be $(b-c)$. So for any
spin grading $s>b-c$, the disks described by the above domains (i.e. between $A_{2b-j,k}$ and $B_{2a-j-1,j+k+\epsilon(j)}$ ) are the only
disks with the described property that are relevant to the boundary map. We will be interested in computing the number of holomorphic
disks whose associated domain is as above.\\

Fix a spin grading $s=b+c-p$. Then $s(A_{2b-j,k})=s$ iff $j+k=p$. These are $A_{2b-1,p-1},A_{2b-2,p-2},....,A_{2b-p,0}$. Correspondingly
we will have the B-type Kauffman states $B_{2a-2,p+1}, B_{2a-3,p},B_{2a-4,p+1},B_{2a-5,p},....$. The disk between $A_{2b-l,p-l}$ and
$B_{2a-1-l,p+\epsilon(l)}$ will be as is shown in \textbf{figure.6}, where $r=2c-p-\epsilon(l)$. We will denote  such a domain by
$\mathcal{D}(l,r)$.\\

\begin{figure}
   \mbox{\vbox{\epsfbox{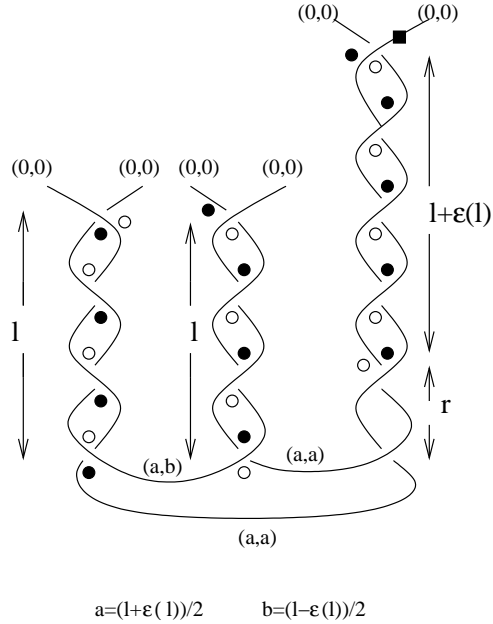}}}
   \caption{\label{fig:figure1}
   {The domain illustrated above is $\mathcal{D}(l,r)$ In order to relate the contribution of this domain to that of $\mathcal{D}(l,r-1)$
    one may consider the homology of the surface obtained by adding long cylinders to the ends of strands and add a $\alpha$-circle
    at the point denoted by a bold square.}}
\end{figure}

\textbf{Lemma.1.} \emph{If $\phi(l,r)$ is the disk associated to $\mathcal{D}(l,r)$ and $\mathcal{M}(l,r)=\mathcal{M}(\phi(l,r))$ is the moduli
space of holomorphic representatives of the disk $\phi(l,r)$, then:
\begin{equation}
\#\{\mathcal{M}(l,r)/\mathbb{R}\}=\pm \#\{\mathcal{M}(l,r-1)/\mathbb{R}\}
\end{equation}}
\\

\textbf{Proof.} In the Heegaard diagram which we made by thickening the knot projection, change the place of the exceptional curve to the point
marked by the solid square in \textbf{figure.6}. Puncture all of the regions which are not shown in the picture. Standard arguments, as those
used by Ozsv\'{a}th and Szab\'{o} for defining the $3$-manifold and Knot Floer homologies show that still we will have a Floer homology, meaning that
there is a boundary operator $\partial$ with the property that $\partial \circ \partial=0$ (one may connect the ends cyclically to get a link /knot
diagram and puncture all of the non-relevant domains. Then put the special $\alpha$-curve at the point marked by a solid square in the picture
and use the results of Ozsv\'{a}th and Szab\'{o}.).\\
The two Kauffman states which are shown in the picture, may be extended to Kauffman states of this diagram, which we will call \textbf{x,y}, with
$s(\textbf{x})=s(\textbf{y})$ and $m(\textbf{x})=m(\textbf{y})+1$. In the same picture we may consider the two Kauffman states which correspond
to $\mathcal{D}(l,r-1)$ and call them \textbf{x$'$,y$'$} with the Maslov index of \textbf{x$'$} higher than that of \textbf{y$'$}. Note that
$s(\textbf{x})-s(\textbf{x}')=1,m(\textbf{x})-m(\textbf{x}')=1$ and that the domain associated to the topological disk connecting these two Kauffman
states is a square, with a circle removed as shown in \textbf{figure.7}. The $\alpha$-circles are shown by solid lines while the $\beta$-circles are the
dashed lines.

It is not hard to see that this domain is associated to a disk, supporting a unique holomorphic representative, as is shown in \cite{OS-alt}.
So by allowing the disks which do not preserve the spin-grading, we will get \textbf{x}$'$ in the boundary of \textbf{x} with coefficient $\pm 1$.
A similar argument shows that  \textbf{y}$'$ appears in the boundary of \textbf{y} with coefficient $\pm 1$.\\

\begin{figure}
   \mbox{\vbox{\epsfbox{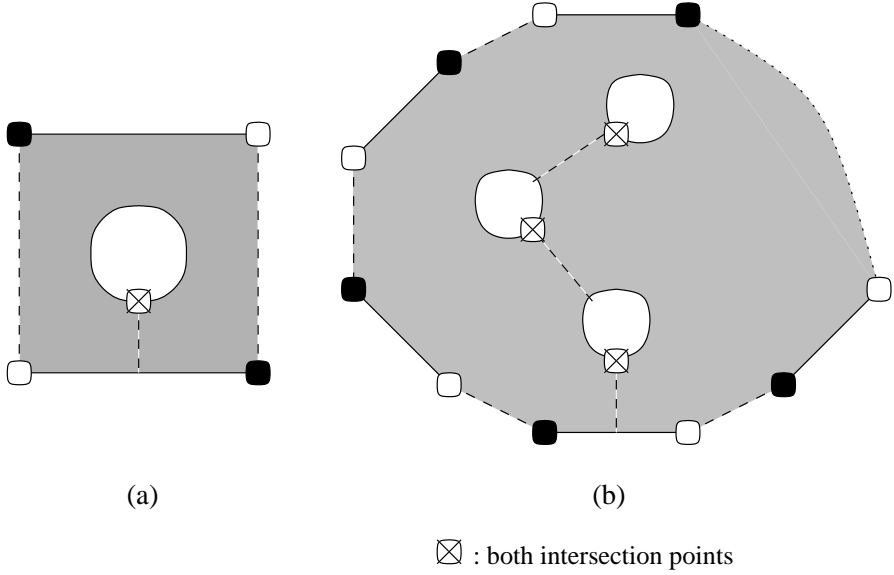}}}
   \caption{\label{fig:figure1}
   {The contribution of the disks with a domain as above, to the boundary map is equal to $\pm 1$}}
\end{figure}

On the other hand, note that $\partial(\textbf{x})=\partial_0(\textbf{x})+\partial_1(\textbf{x})+\partial_2(\textbf{x})+....$ where
$\partial_i(\textbf{x})$ is the part of boundary of \textbf{x} whose spin grading is equal to that of \textbf{x} minus $i$. Since
$\partial(\partial(\textbf{x}))=0$ the coefficient of \textbf{y}$'$ in this expression, denoted by $c_{\textbf{y}'}(\partial(\partial(\textbf{x})))$
is also $0$. Consequently the above coefficient in $\partial(\partial_0(\textbf{x})+\partial_1(\textbf{x}))$ is
zero, since the boundary operator does not increase the spin grading. It is easy to see that the only relevant Kauffman states to
$\partial_0(\textbf{x})$ and $\partial_1(\textbf{x})$ are respectively \textbf{y,x}$'$. Let $\partial_0(\textbf{x})=\lambda \textbf{y}$. Then
\begin{equation}
0=\partial (\partial (\textbf{x}))=c_{\textbf{y}'}(\partial(\partial(\textbf{x})))=c_{\textbf{y}'}(\partial (\lambda \textbf{y}\pm \textbf{x}'))=
\pm \lambda \pm c_{\textbf{y}'}(\partial(\textbf{x}'))
\end{equation}
This completes the proof of the lemma.$\square$\\

\section{Main Theorems}
Our goal is to prove the following theorem from which we will deduce our claim about the homology of pretzel knots:\\

\textbf{Theorem.3.} \emph{With the above notation, the contribution of the domain $\mathcal{D}(l,r)$ to the boundary map is $\pm 1$. Saying it differently:
\begin{equation}
\#\{\mathcal{M}(l,r)/\mathbb{R}\}=\pm 1
\end{equation}}
\\

Fist of all note that this theorem implies that the pairs of Kauffman states of type (b) in the spin gradings $b-c< s\leq b+c+1$ are canceled
against each other. This means that if there are $m$ Kauffman states of type $A$ (i.e $A_{jk}(b)$), and $n$ of type $B$, then there will be
a group $\oplus_{i=1}^{m-n} \mathbb{Z}$ of Maslov index $m=s-b-c$ or a group $\oplus_{i=1}^{n-m} \mathbb{Z}$ of Maslov index $m=s-b-c-1$ if
$n\leq m$ or $m\leq n$ respectively. This determines the groups for $s$ in the above range. We will comment on the case $s\leq b-c$ later.\\

\begin{figure}
   \mbox{\vbox{\epsfbox{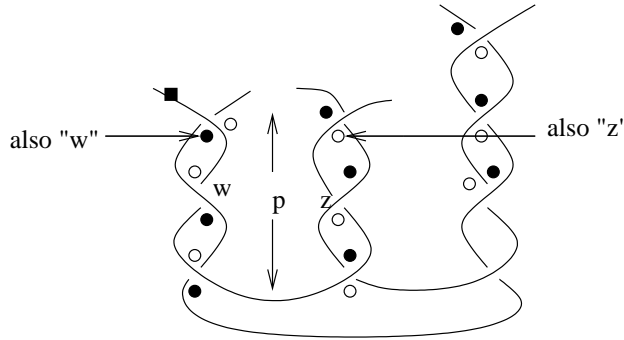}}}
   \caption{\label{fig:figure1}
   {The contribution of disks with domains of type $\mathcal{D}(p,r)$ for $p$ even, may be computed directly from that of
    $\mathcal{D}(p-1,r)$}}
\end{figure}

\textbf{Proof of theorem.3.} Let us put $S(l)$ for the statement that the domain $\mathcal{D}(l,r)$ supports a unique (up to reparametrization)
holomorphic disk. Lemma.1 says that it does not matter which $r$ we are considering and in particular we may restrict ourselves to the case $r=0$.
We will show $S(l)$ by induction on $l$. For $l=0$, this is our rectangle, and there is nothing to prove. Now suppose that
$S(0),...,S(p-1)$ are true. Consider the spin-grade $s=b+c-p$ in the type (a) Heegaard diagram. The relevant A, B and C type Kauffman states are
$A_{2b-j,p-j}(a), j=1,2,...,p$ (note that $p<2b$), $B_{(2a-1)-j,p+\epsilon(j)}(a), j=1,2,...,2a-2$ and $C_{0,2b-(p+1)}(a)$. There is a disk between
$A_{2b-p,0}(a)$ and $C_{0,2b-(p+1)}(a)$, which is of the form described in \textbf{figure.7(b)}. In fact it will be an $8$-gon with $a+b-p$ holes.
The disks
of the type shown in \textbf{figure.7} are easily seen to have a unique holomorphic representative (This may be done by considering (-2,3,2a+1)-pretzel
knots, where the computation may be done with an indirect approach, at least for the $8$-gon). Although it uses two of the regions on the top, these are
on the upper face of the surface, while for the above discussions, we only  punctured the lower faces in order to get rid of the rest of the
disks. As a result, we may assume that $A_{2b-p,0}(a)$ and $C_{0,2b-(p+1)}(a)$ are canceled against each other. The disks between the pairs
$(A_{2b-j,p-j}(a),B_{(2a-1)-j,p+\epsilon(j)}(a))$ are of the form $\mathcal{D}(j,r)$ for some number $r$ and $j=1,...,p-1$, so from the induction
hypothesis, we know that
these pairs may also be assumed to cancel each other. If we have been careful to choose $a$ such that $2a=p$ or $2a-1=p$ (Note that for the proof
of this theorem we are free to change our numbers $a,b,c$) then all of the B-type Kauffman states will be canceled and the homology groups will be
completely supported in degree $-p$. Note that there is no need to consider disks with positive coefficients in the punctured regions since there is no
more possible cancellation.\\

Now we do the same computation, but with the (b)-type diagram. In this case the relevant Kauffman states are
$A_{2b-j,p-j}(b), j=1,...,p$  and $B_{(2a-1)-j,p+\epsilon(j)}(b), j=1,2,...,2a-1$. Again all of B-type states are canceled against respective A-type
states because of induction hypothesis, except possibly for $B_{0,p+1}(b)$. Consider the following two possible cases:\\
\\
(a) $p$ is odd. Then $a$ may be chosen such that $p=2a-1$. Then after the above cancellation the only Kauffman states that remain are $A_{2b-p,0}(b)$ and
$B_{0,p+1}(b)$. The connecting disk is of the type $\mathcal{D}(p,r)$ and from the above computation of the Floer homology we know that the Floer homology
is zero and therefore there should be a unique holomorphic representative and we are done.\\
\\
(b) $p$ is even. In this case we will prove the theorem by a direct argument. The proof will be very similar to that of \textbf{lemma.1}.
We will fix a reference point as shown by a bold square in \textbf{figure.8} and consider the relevant Floer homology, allowing the
disks that change spin-grading. Again we may introduce two new Kauffman states \textbf{w,z}, which agree with our white and bold states
respectively, in all but two intersections, where the difference is illustrated in \textbf{figure.8}. For convenience we will call the
white and bold state \textbf{y,x} respectively. $s(\textbf{w})=s(\textbf{z})=s(\textbf{x})+1$ and again there are disks with unique
holomorphic representatives going from \textbf{w} to \textbf{y} and from \textbf{z} to \textbf{x}. By considering the coefficient of
\textbf{y} in $\partial \circ \partial (\textbf{z})$ we see that:
\begin{equation}
c_{\textbf{y}}(\partial \textbf{x})=\pm c_\textbf{w}(\partial \textbf{z})
\end{equation}
But the disk between \textbf{z} and \textbf{w} is of type $\mathcal{D}(p-1,r)$ and we are done by the induction hypothesis. Note that
it is important here that $p$ is even and the same argument will not work for odd $p$. $\square$   \\
\\

\textbf{Proof of theorem.1:} Note that the above argument will give the parts 1,2,3 of the theorem as noted before the proof of the previous theorem.
Parts 5,6 will follow because of the symmetry. In fact there will be exactly similar disks between $A$ and $C$-type Kauffman states. For the part
4, there are two different cases, one for $p=2b, ..., 2c-1$ and the other one for $p=2c$. Again for the first case there are only A and B type Kauffman
states with any B type state being paired with an A type state and the disk between these two being of the type described in lemma.1 and theorem.3.
So they are canceled against each other. For the second case, we will have A,B and C type Kauffman states. But still, all of the B and C type ones
have associated A type states which cancel them. This completes the proof of theorem.1. $\square$\\
\\
The proof of theorem.2 is essentially similar. Also one may state a theorem like theorem.2 for $P(-2a,2b+1,2c+1)$ when $a$ is not necessarily less than
$b,c$. The essential tool in all proofs will be lemma.1 and theorem.3.


\begin{thebibliography}{99}
\bibitem{knot} Likorish, W. B. R., \emph{An Introduction to knot theory}, Graduate texts in mathematics, 175, Springer-Verlag,
New York, 1997
\bibitem{OS-knot}  Ozsv\'{a}th, P. and Szab\'o, Z., \emph {Holomorphic disks and knot invariants}, preprint, available at math.GT/0209056
\bibitem{OS-3m1}  Ozsv\'{a}th, P. and Szab\'o, Z., \emph {Holomorphic disks and topological invariants for closed three-manifolds}, to
appear in Annals of Math.,
available at math.SG/0101206
\bibitem{OS-3m2}  Ozsv\'{a}th, P. and Szab\'o, Z., \emph {Holomorphic disks and three-manifold invariants: properties and applications},
to appear in Annals of Math.,
available at math.SG/0105202
\bibitem{OS-4m} Ozsv\'{a}th, P. and Szab\'o, Z., \emph {Holomorphic triangles and invariants for smooth four-manifolds} preprint,
available at math.SG/0110169
\bibitem{OS-fibered} Ozsv\'{a}th, P. and Szab\'o, Z., \emph {Heegaard Floer homologies and contact structures}, preprint,
available at math.SG/0210127

\bibitem{OS-alt}  Ozsv\'{a}th, P. and Szab\'o, Z., \emph {Heegard Floer homology and alternating knots}, preprint, available at math.GT/0209149
\bibitem{Ras} Rasmussen, J., \emph{Floer homology of surgeries on two bridge knots}, Alg. and Geo. Topology,2 (2002) 757-789
\bibitem{Ras2} Rasmussen, J., \emph{Floer homology and knot complements}, Ph.D thesis,
also available at math.GT/0306378
\end{thebibliography}
\end{document}